\documentclass{amsart}

\usepackage{graphicx}

\begin{document}

\title[Tile Count of Regular $2n$-gons]{Tile Count in the Interior of Regular $2n$-gons Dissected by Diagonals Parallel to Sides}

\author{Richard J. Mathar}
\urladdr{http://www.strw.leidenuniv.nl/~mathar}
\email{mathar@strw.leidenuniv.nl}
\address{Leiden Observatory, Leiden University, P.O. Box 9513, 2300 RA Leiden, The Netherlands}
\thanks{This work was inspired by sequence A165217 of the OEIS\cite{EIS}.
}


\date{\today}

\begin{abstract}
The regular $2n$-gon (square, hexagon, octagon,\ldots) is subdivided into
smaller polygons (tiles) by the subset of diagonals which run parallel to any of the
$2n$ sides. The manuscript reports on the number of tiles up to the 78-gon.

\end{abstract}

\maketitle
\section{Summary} 
Given the $N=2n$ sided regular polygon, its interior is dissected into
non-overlapping regions (polygons, tiles) by the $n(n-2)$ diagonals
parallel to any of the polygon's sides that connect edges of the polygon.
The number of unrestricted diagonals is $n(2n-3)$. $n(n-1)$ of these are not
considered here, so the number of regions $F$ remains smaller than those generated
by all of them \cite{PoonenDM11,HarborthEM24,SchneiderEM14}.

The $N$ sides define $N/2$ different non-oriented directions, each represented
by a bundle of $N/2-2$ parallel diagonals; the subtraction of 2 indicates
the polygon sides are not included in the count.
The tile counts are summarized in Tables \ref{tab.c28} and \ref{tab.c30}.

\begin{table}[h]
\begin{tabular}{lrrrrrrrrrrrrr}
$N$ &   4 & 6 & 8 & 10 & 12 & 14 & 16 & 18 & 20 & 22 & 24 & 26 & 28 \\
$n$ &   2 & 3 & 4 &  5 &  6 & 7  &  8 &  9 & 10 & 11 & 12 & 13 & 14 \\
\hline
tiles $F$&  1 & 6 & 25& 50& 145 & 224&497 & 630 &1281& 1606 & 2761 & 3302 & 5265\\
edges $E$&  4 & 12 & 48& 80& 276& 378& 960 & 1062 &2500& 2860 & 5424 & 5980& 10388\\
\end{tabular}
\caption{Common results of the manual count of Figures \ref{fig6.ps}--\ref{fig28.ps} and
of the \texttt{C++} program of the appendix.
}
\label{tab.c28}
\end{table}

\begin{table}[h]
\begin{tabular}{lrrrrrrrrrr}
$N$ &   30 & 32 & 34 & 36 & 38 & 40 & 42 & 44 & 46 & 48 \\
$n$ &   15 & 16 & 17 & 18 & 19 & 20 & 21 & 22 & 23 & 24 \\
\hline
$F$&5940 & 9185 &10472 &14977 &16834 &23161 &25284 &34321& 37720 &49105\\
$E$&10770 &18176&19482 &29700 &31616 &46000 &47460 &68244 &71714 &97728\\
\end{tabular}
\begin{tabular}{lrrrrrrrrr}
$N$ &   50 & 52 & 54 & 56 & 58 &60 &62 &64 &66\\
$n$ &   25 & 26 &27 & 28 &29 &30 &31 &32 & 33\\
\hline
$F$& 53500&68225 & 73278 & 92457 & 99470 &122641 &131316 &159681 &169158\\
$E$&102150&135876 &140076 &184240 & 191284 &244500 &253270 &318464 &326238\\
\end{tabular}
\begin{tabular}{lrrrrrrrrr}
$N$ &   68 & 70 & 72 &74 &76 &78\\
$n$ &   34 & 35 & 36 &37 &38 &39\\
\hline
$F$& 204545 &217210 & 258265 &273282 &321937 &338208\\
$E$& 408068 &420840 &515376 &530432 &642580 &656526\\
\end{tabular}
\caption{Results of the \texttt{C++} program complementing Table \ref{tab.c28}.
}
\label{tab.c30}
\end{table}

By the discrete rotational symmetry,
the polygon is invariant if rotated by multiples of the angle $2\pi/N$. One can count the
nonequivalent tiles in a segment in one particular ray direction,
then multiply by $N$, and increase the count by 1
to include the central tile if $n$ is even.

\section{Illustrations} 
To verify the counts of Table \ref{tab.c28}, a visual inspection follows
for the cases up to the 28-gon. The tiles are enumerated in two different colors
to demonstrate that each is in exactly one of the $2n$ replicas.
Each of the numbers 1, 2, 3, $\ldots, \lfloor F/N\rfloor$
appears $N$ times in a figure. From Figure \ref{fig24.ps} on, the enumeration
covers only three rays to simplify the inspection of
crossing patterns in the cluttered areas.

The standard coordinate system puts the polygon corners on the unit
circle at Cartesian
coordinates
$(\cos(\pi k/n),\sin(\pi k/n))$ for $k=0,\ldots,N-1$. If these can
be expressed as square roots \cite{ConwayDiCompG22,GirstmairAA81,GurakMathComp75,ServiAMM110,BrackenIJQC90,VidunasKJM59},
the coordinates of the diagonal crossings are also of this form, since only determinental
combinations of the coordinates of the polygon corners are involved to determine their positions.

\begin{figure}[hbt]
\includegraphics[scale=0.5,clip]{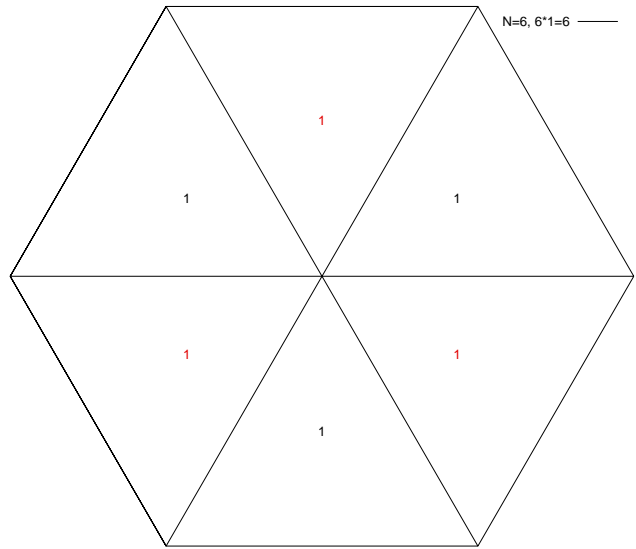}
\caption{
$N=6$ sides: 6 tiles, one triangular tile replicated
6 times. Corner coordinates are $(\pm 1,0)$, $(\pm 1/2,\pm \sqrt{3}/2)$.
}
\label{fig6.ps}
\end{figure}

\begin{figure}[hbt]
\includegraphics[scale=0.6,clip]{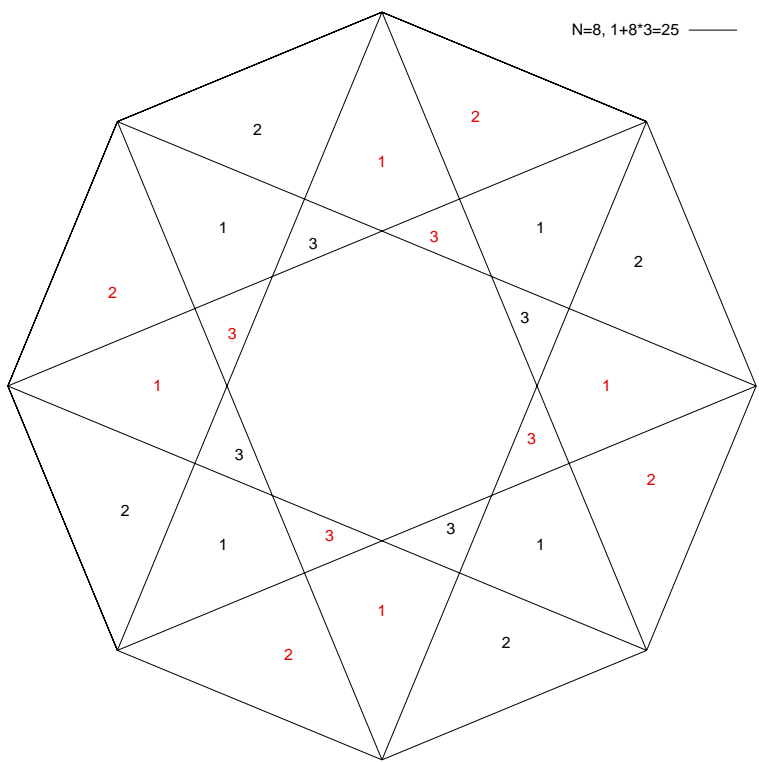}
\caption{
$N=8$ sides: 25 tiles,  one octagonal tile in the center and 3 triangular
or kite-shaped tiles replicated 8 times outside the center.
The 8 corner coordinates on the unit circle
are $(\pm 1,0)$, $(0,\pm 1)$, $(\pm \sqrt{2}/2,\pm \sqrt{2}/2)$. Internal crossings occur at $(\pm 1/2, \pm (\sqrt{2}-1)/2)$,
$(\pm (1-\sqrt{2}),0)$, $(0,\pm(1-\sqrt{2}))$, $(\pm(1-\sqrt{2}/2),\pm(1-\sqrt{2}/2)$ and $(\pm(1-\sqrt{2})/2,\pm 1/2)$.
}
\end{figure}

\begin{figure}[hbt]
\includegraphics[scale=0.8,clip]{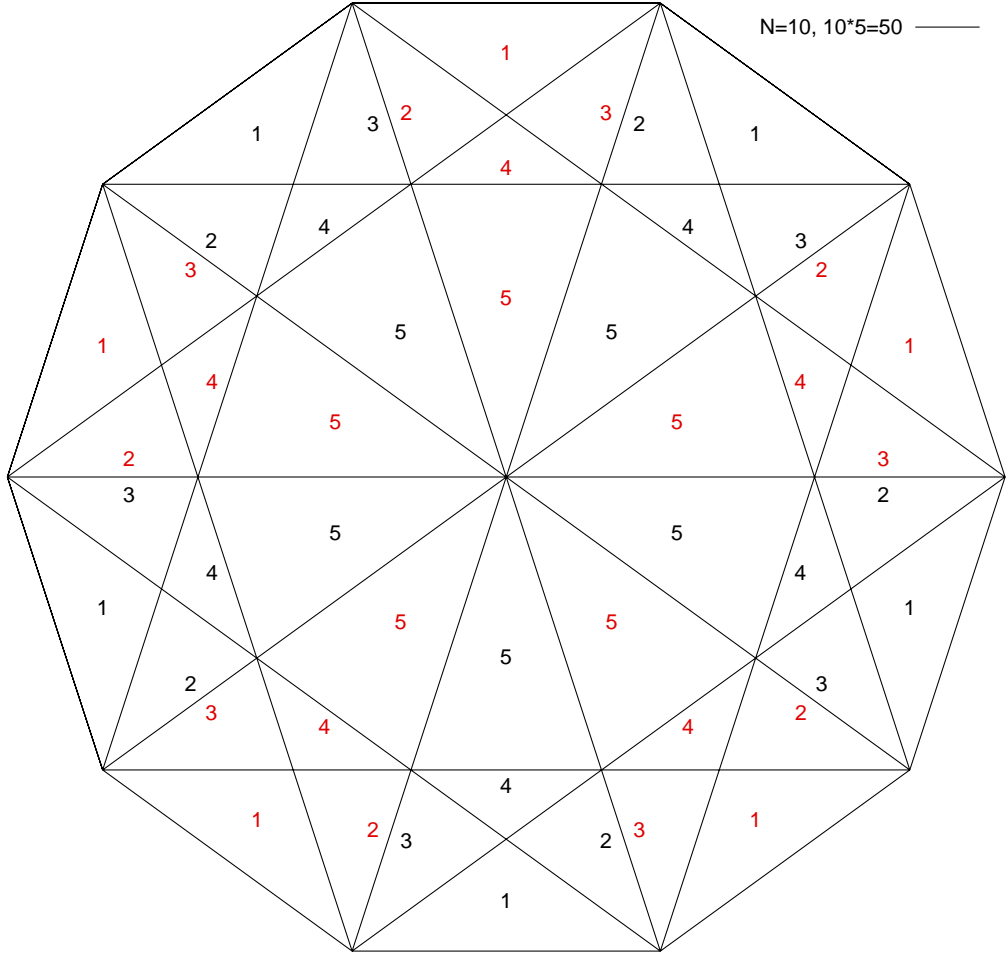}
\caption{
$N=10$ sides: 50 tiles, comprising 5 triangular tiles replicated
10 times.
The 10 corner coordinates on the unit circle are
$(\pm 1,0)$,
$(\pm (1+\sqrt{5})/4,\pm \sqrt{2(5-\surd{5})}/4)$,
$(\pm (\sqrt{5}-1)/4,\pm \sqrt{2(5+\surd{5})}/4)$.
Internal crossings are at
$(\pm (3-\sqrt{5})/4,\pm \sqrt{2}\sqrt{5-\surd{5}}/4)$,
$(\pm 1/2,\pm \sqrt{5-2\surd{5})})$,
$(\pm (5-\sqrt{5})/4,\pm (\sqrt{5}-2)\sqrt{2(5+\surd{5})}/4)$,
$(\pm (1+\sqrt{5})/4,\pm \sqrt{2(5-\surd{5})}/4)$,
$(0,\pm (3-\sqrt{5})\sqrt{2(5-\surd{5})}/4)$,
$(\pm (1-\sqrt{5})/2,0)$ and $(0,0)$.
}
\end{figure}

\begin{figure}[hbt]
\includegraphics[scale=0.8,clip]{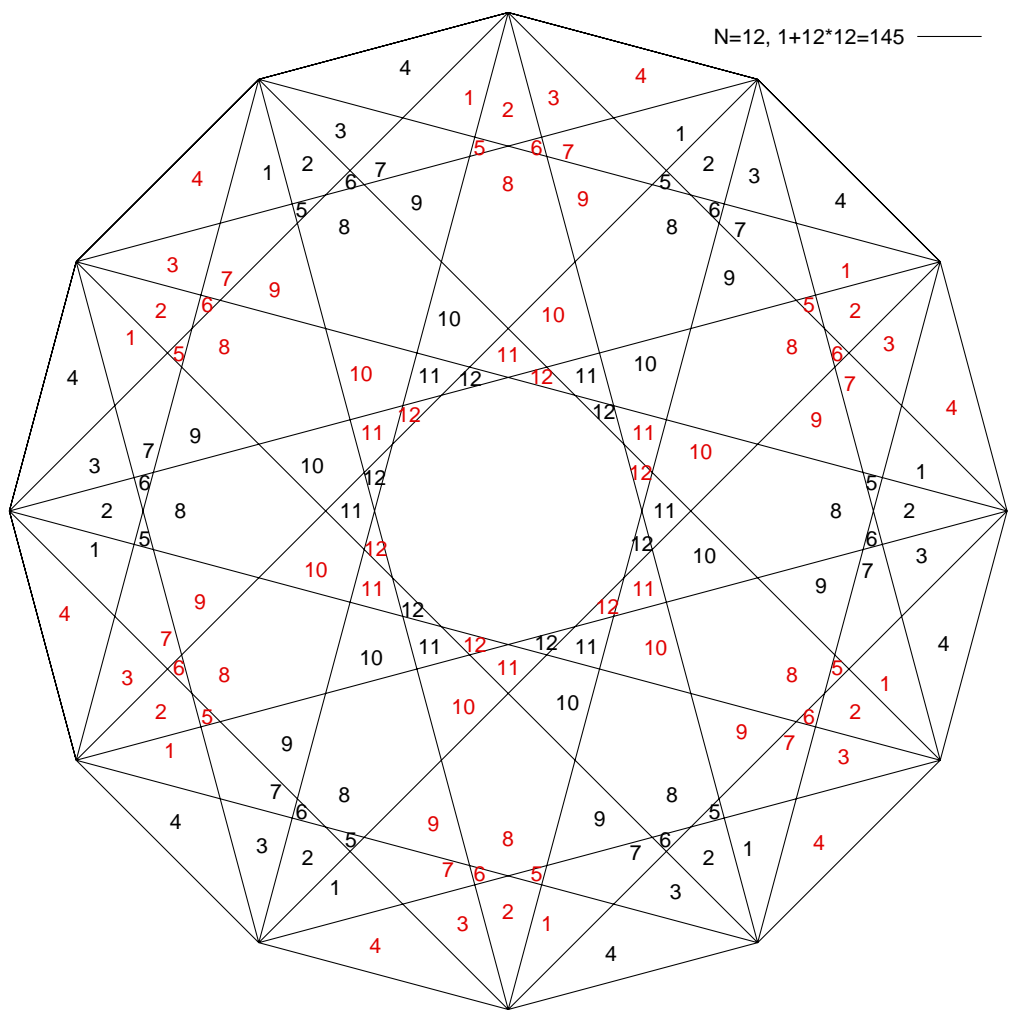}
\caption{
$N=12$ sides: 145 tiles. One tile in the center and 12 tiles replicated
12 times.
}
\end{figure}

\begin{figure}[hbt]
\includegraphics[scale=0.8,clip]{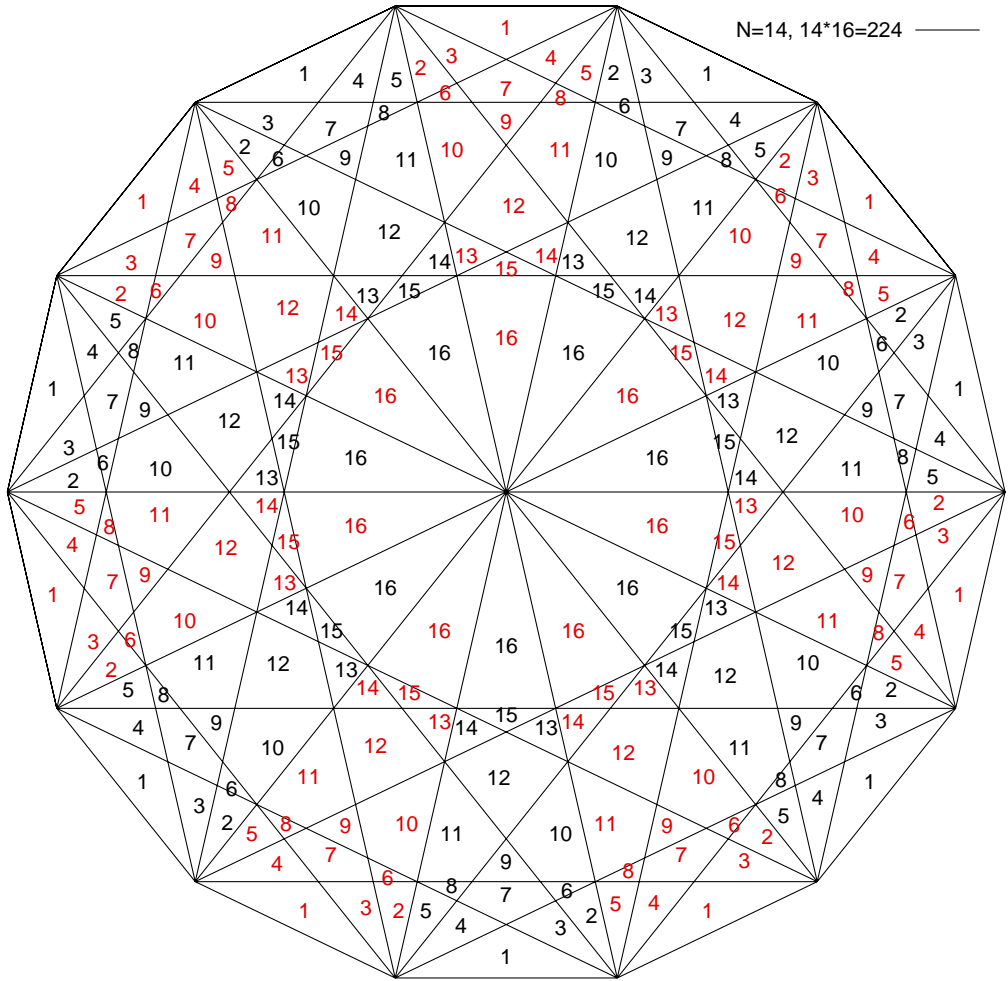}
\caption{
$N=14$ sides: 224 tiles, comprising 16 tiles replicated
14 times.
}
\end{figure}

\begin{figure}[hbt]
\includegraphics[scale=0.95,clip]{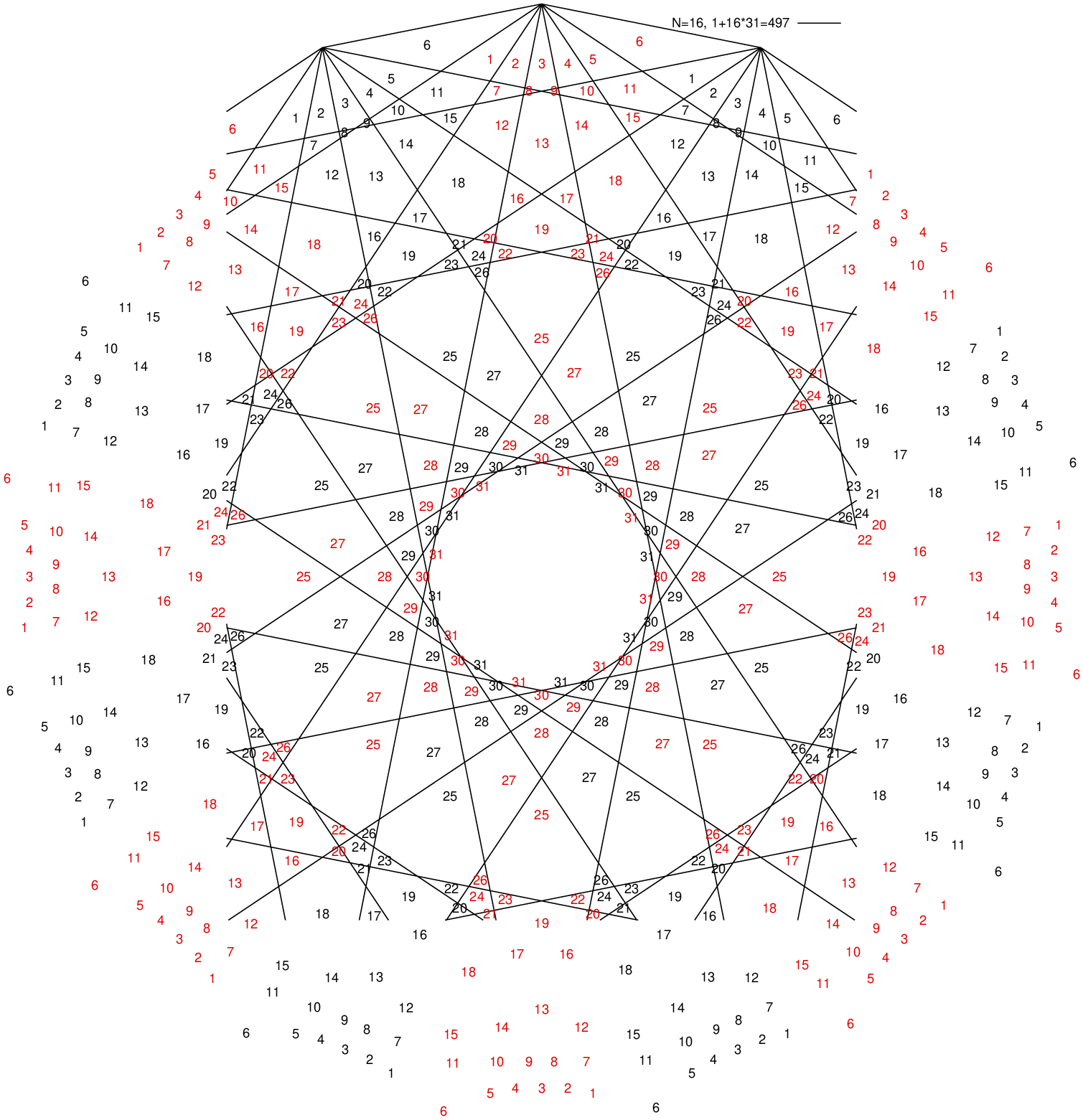}
\caption{
$N=16$ sides: 497 tiles. One tile in the center and 31 tiles replicated
16 times after rotation.
}
\end{figure}

\begin{figure}[hbt]
\includegraphics[scale=0.95,clip]{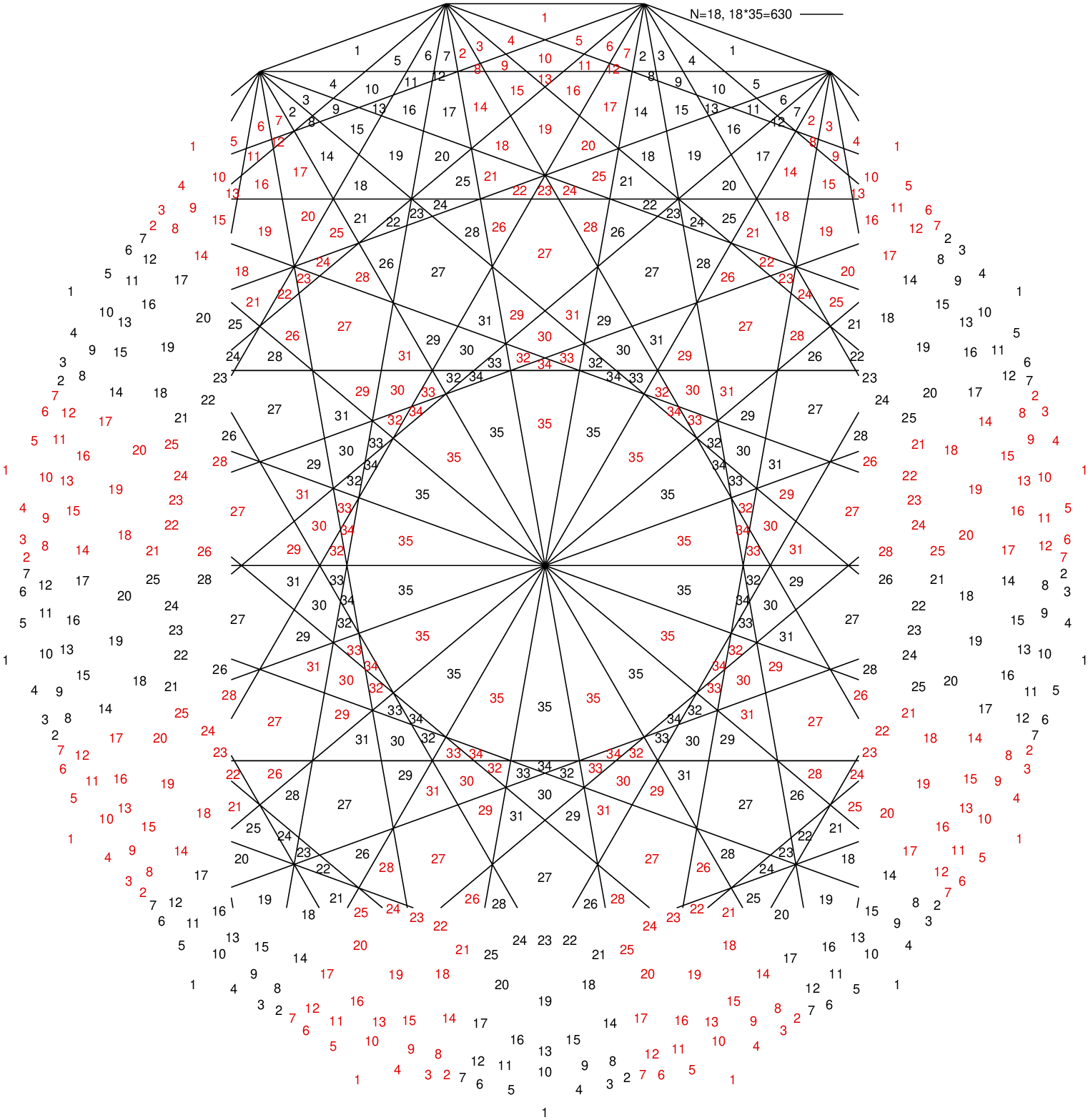}
\caption{
$N=18$ sides: 630 tiles, consisting of 35 tiles replicated
18 times.
}
\end{figure}

\begin{figure}[hbt]
\includegraphics[scale=0.9,clip]{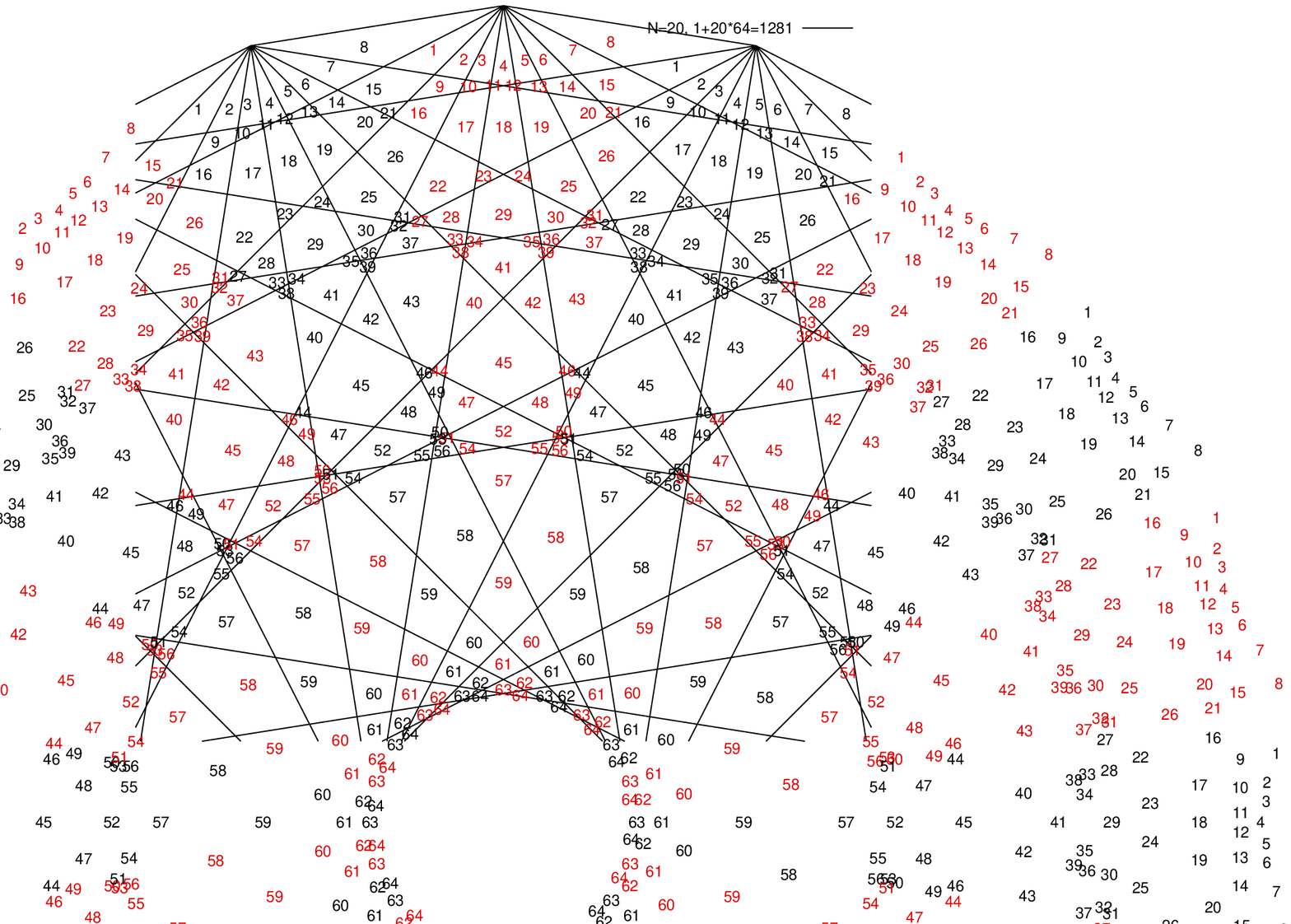}
\includegraphics[scale=0.9,clip]{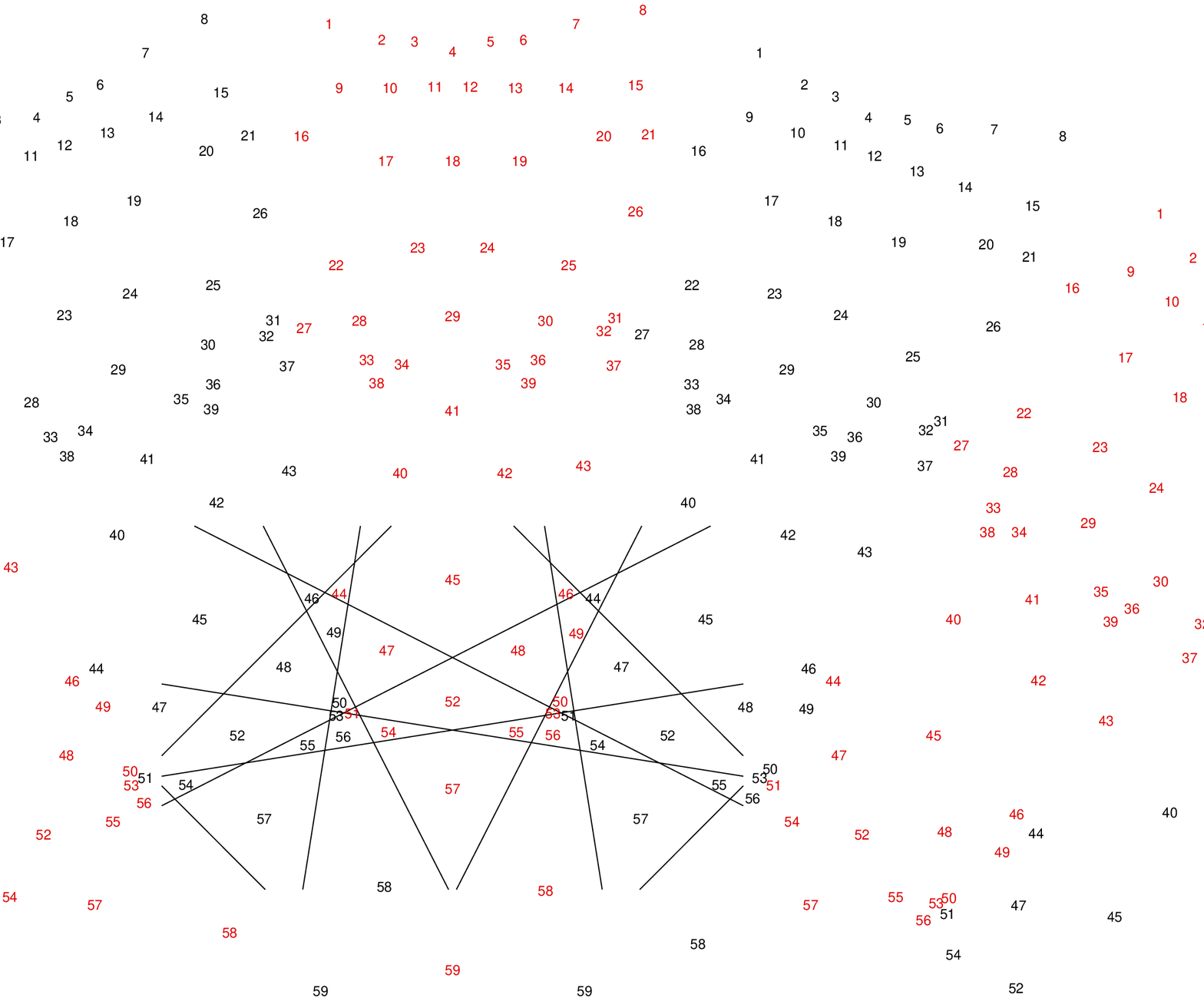}
\caption{
$N=20$ sides: 1281 tiles. One tile in the center and 64 tiles replicated
20 times.
The lower graph shows a zoomed section of the upper graph.
}
\label{fig20.ps}
\end{figure}

\begin{figure}[hbt]
\includegraphics[scale=0.95,clip]{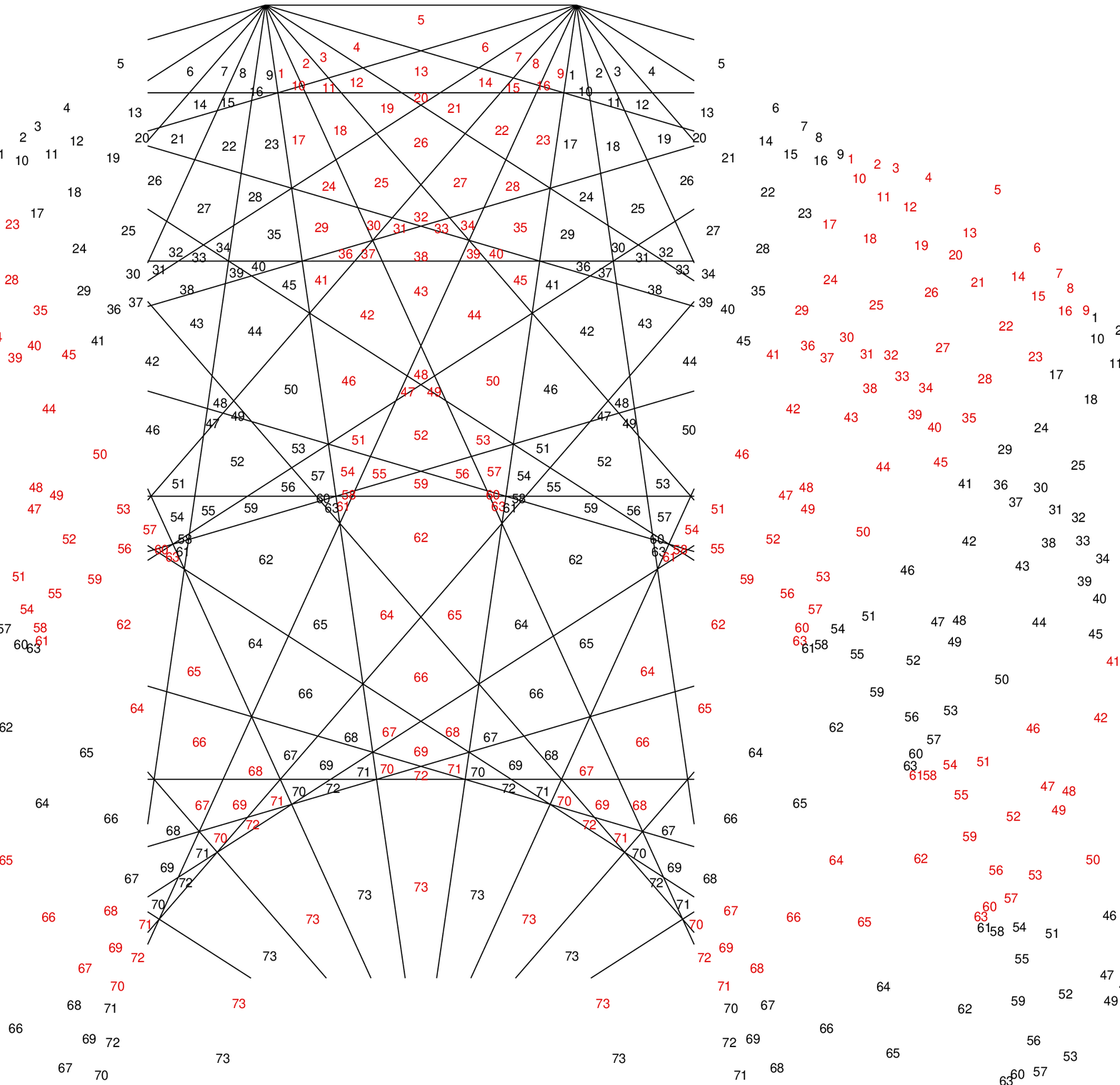}
\caption{
$N=22$: 1606 tiles, consisting of 73 tiles replicated
22 times.
}
\label{fig22.ps}
\end{figure}

\begin{figure}[hbt]
\includegraphics[scale=0.95,clip]{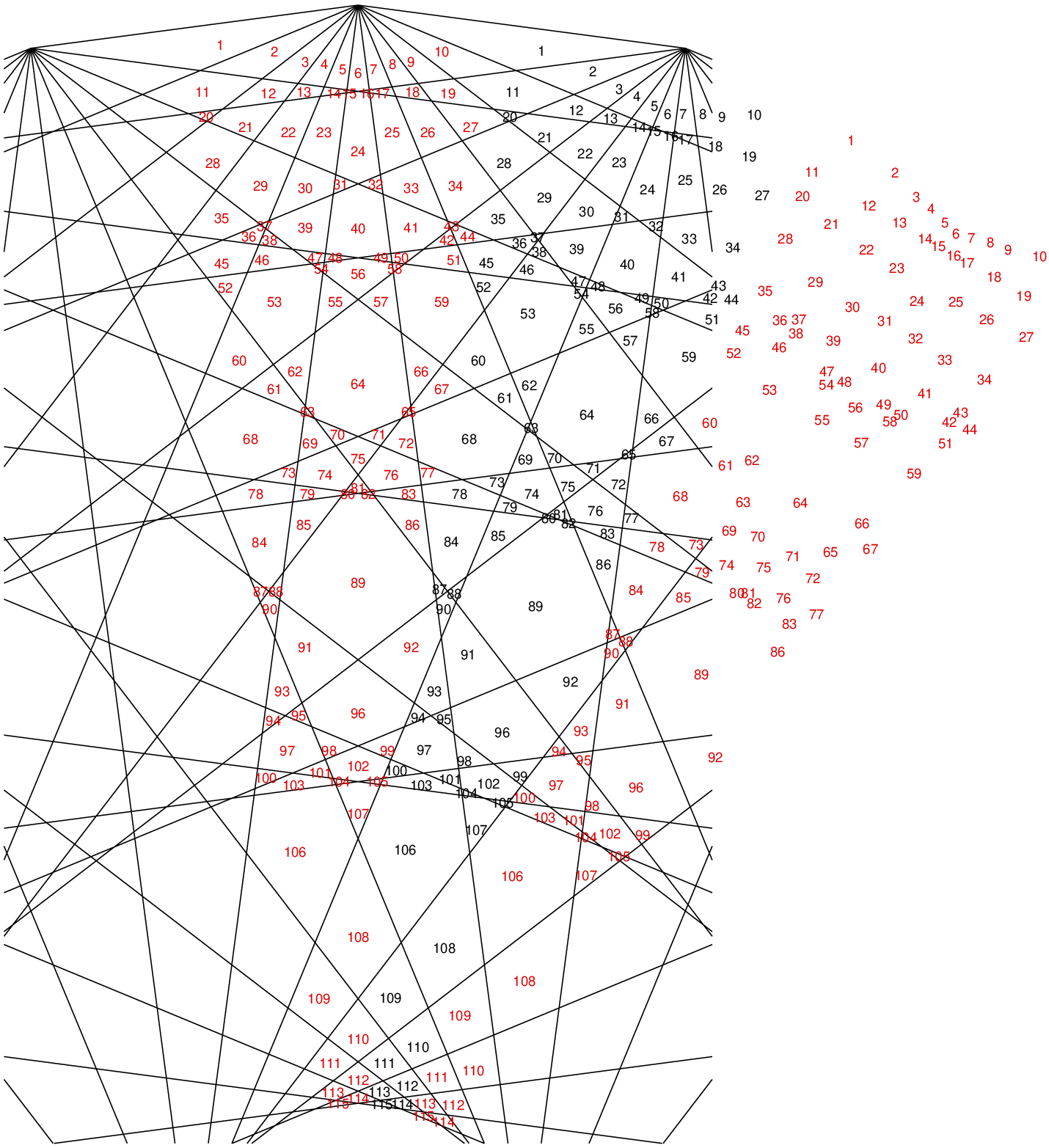}
\caption{
$N=24$: 2761 tiles, one in the center and 115 tiles  replicated
24 times after incremental rotation.
}
\label{fig24.ps}
\end{figure}

\begin{figure}[hbt]
\includegraphics[scale=0.97,clip]{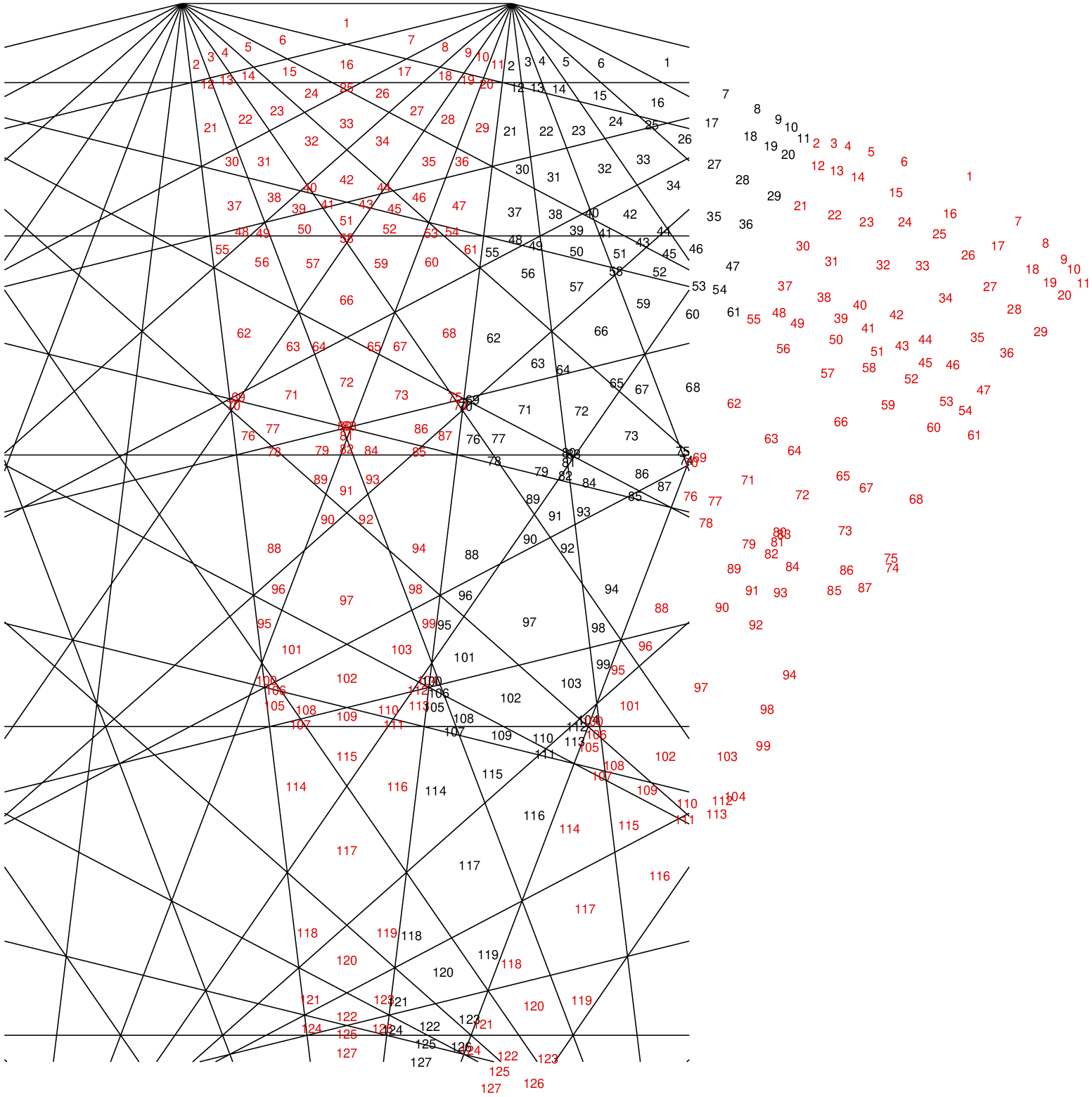}
\caption{
$N=26$:  3302 tiles, 127 tiles per ray replicated
$N$ times.
}
\end{figure}

\begin{figure}[hbt]
\includegraphics[scale=0.9,clip]{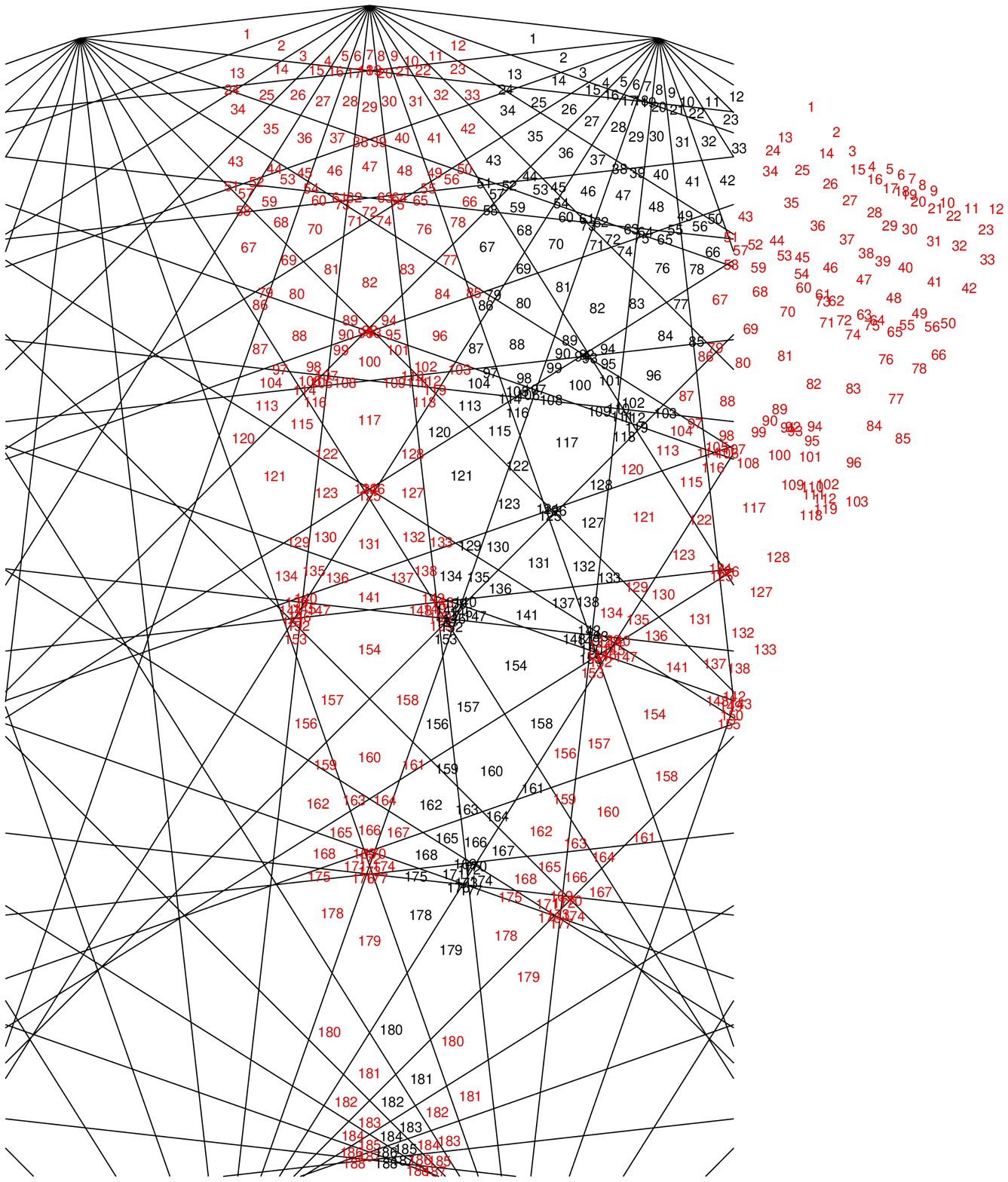}
\caption{
$N=28$: 5256 tiles, one in the center and 188 replicated $N$ times.
}
\label{fig28.ps}
\end{figure}

\clearpage
\appendix
\section{C++ Line Intersection Program}

Another way of counting is to enumerate all $E$ edges and all
$V$ intersections (vertices) and to use the planar version
of Euler's formula
\begin{equation}
F=1+E - V
\end{equation}
to count the faces of the graph. (The presence of the $1$ indicates that
the area outside the polygon is not counted.) The program starts from the edges of
the exterior lines of the polygon that connect the roots of unity in the
complex plane, plus the edges of the diagonals parallel to any of these
edges.  This base set is then rescanned producing a working set of chopped
edges by looping over the segments of the base set,
\begin{enumerate}
\item
each time computing all
points of intersection (including ``touches'') with any element  of the working set,
\item
if this splits the element of the working set, replacing the element in
the working set by its fragments, 
\item
adding the base set element or (if intersected or touched) its sub-components to 
the working set.
\end{enumerate}
The loop finished, $E$ becomes the number of elements in the working set.
$V$ is computed by gathering all $2E$ terminal points of these in a
point set, discarding duplicates that have a mutual distance smaller
than some noise threshold set by the floating point arithmetic.

The program is simple; no attempt is made to take advantage of the $N$-fold
rotational symmetry of the graph.

\scriptsize \begin{verbatim}
#include <iostream>
#include <cmath>
#include <vector>
#include <set>
#include <algorithm>

using namespace std ;

/** A point with two Cartesian coordinates in the plane.
*/
class Point {
public:
    /** Abscissa doordinate.
    */
    double x ;
    /** Ordinate coordinate.
    */
    double y ;

    /** Default ctor. Point at infinity.
    */
    Point() : x(INFINITY), y(INFINITY)
    {
    }

    /** Ctor from two known Cartesian coordinates.
    * @param xcoo The x coordinate.
    * @param ycoo The y coordinate.
    */
    Point(const double xcoo, const double ycoo) : x(xcoo), y(ycoo)
    {
    }

    /** Distance to another point.
    * @param oth The other point to compare with.
    * @return The Euclidean distance between this point and oth.
    */
    double dist(const Point & oth) const
    {
        return hypot( x-oth.x, y-oth.y) ;
    }

    /** A distance somewhat larger than the estimated floating point error in the
    * computation of point positions.
    */
    static const double pointFuzzy = 1.e-10 ;
protected:
private:

} ; /* Point */

/** Comparison of two points.
* Two points are considered equal if their mutual distance is (in the current coordinate units)
* smaller than the constant value defined by Point::pointFuzzy .
*/
bool operator== (const Point &a, const Point b)
{
    return ( a.dist(b) <= Point::pointFuzzy );
}

/** Comparison of two points. 
* This operator supports operations on point sets of the STL, and the underlying ordering
* is not relevant for this program here.
*/
bool operator< (const Point &a, const Point & b)
{
    if ( a.dist(b) <= Point::pointFuzzy )
        return false;
    else if ( a.x < b.x-Point::pointFuzzy)
        return true;
    else if ( a.x > b.x+Point::pointFuzzy)
        return false;
    else
        return a.y < b.y ;
}

/** A straight line segment.
* A connection between two points.
*/
class Line {
public:
    /** The two terminal points of the segment.
    */
    Point pts[2] ;

    /** Ctor given a start and end point.
    * @param strt One of the two terminal points.
    * @param fin The terminal points at the opposite end.
    */
    Line(const Point & strt, const Point & fin)
    {
        pts[0] = strt ;
        pts[1] = fin ;
    }

    /** Copy ctor.
    * @param oth The other line that defines the new instance.
    */
    Line( const Line &oth)
    {
        pts[0] = oth.pts[0] ;
        pts[1] = oth.pts[1] ;
    }

    /** Euclidean length.
    * @return Pythagorean distance between the two terminal points.
    */
    double len() const
    {
        return pts[0].dist(pts[1]) ;
    }

    /** Determine the point on the line defined by parameter t.
    * @param t The parameter along the line in units of the segment length. Values of 0 and 1 recall the terminal points.
    * @return The point at pts[0] + t*(pts[1]-pts[0]).
    * If t was outside the interval [0,1], this point is outside the line segment, but still coplanar.
    */
    Point atParam(const double t) const
    {
        return Point( t*pts[1].x +(1.-t)*pts[0].x,t*pts[1].y +(1.-t)*pts[0].y ) ;
    }

    /** Detect intersection with another line.
    * @param[in] oth  The other line to intersect with.
    * @param[out] tOth The point of intersection parametrized by oth.strt+tOth*(oth.fin-oth.strt). 
    * @return the point of intersection given by the parametrization t in strt+t*(fin-strt).
    * If the lines are parallel, a value of infinity is returned. This line segment here intersects with 'oth'
    * if 0<=t<=1 and 0<=tOth <=1. Touching line segments have one of these two parameters in [0,1] and the other
    * one very close to 0 or 1.
    */
    double inters(const Line &oth, double & tOth) const
    {
        /* Computation solves for the two linear equations for the point of intersection, and looks at the
        * determinant to detect parallelism.
        */
        double A[2][2] ;
        double rhs[2] ;
        A[0][0] = pts[1].x - pts[0].x ;
        A[0][1] = oth.pts[0].x - oth.pts[1].x ;
        A[1][0] = pts[1].y - pts[0].y ;
        A[1][1] = oth.pts[0].y - oth.pts[1].y ;
        rhs[0] = oth.pts[0].x - pts[0].x ;
        rhs[1] = oth.pts[0].y - pts[0].y ;

        /* determinant of the coefficient matrix A of the two linear equations */
        const double det = determinant(A[0][0],A[0][1],A[1][0],A[1][1]) ;

        /* Use some measure of the determinant relative to the length
        * of the two vectors to put a threshold in a floating point environment.
        */
        if ( fabs( det) < pointFuzzy*fabs(len()* oth.len()) )
        {
            tOth = INFINITY ;
            return INFINITY ;
        }
        else
        {
            /* If the determinant does not vanish, Cramer's rule computes the two line paremeters t and t0th.
            */
            tOth = determinant(A[0][0],rhs[0],A[1][0],rhs[1])/det ;
            return determinant(rhs[0],A[0][1],rhs[1],A[1][1])/det ;
        }
    }

    
    /** Compute the two line segments if this line is split at parameter t
    * @param t The line parameter in the range (0,1).
    * @return The two line segments. The first one representing the (0,t) interval of points,
    *  the second one representing (t,1).
    */
    vector<Line> splitAt(const double t) const
    {
        vector<Line> out ;
        if ( inSplitRange(t) )
        {
            /* create and attach the two line segments */
            out.push_back( Line(pts[0], atParam(t) ) );
            out.push_back(  Line(atParam(t) ,pts[1] ) );
        }
        else
            /* invalid specification. Return one component only, the segment itself.
            */
            out.push_back(  Line(*this) );
        return out ;
    }

    /** Decompose the line into non-overlapping segments defined by a set of line parameters of the new terminal points.
    * @param t The vector with numbers between 0 and 1 that define the locations of the splitting location.
    * @return The segmented version. The union of all these segments is the line itself.
    * @warn The current implementation does not check that all t[] are in the range [0,1].
    */
    vector<Line> splitAt( vector<double> t) const
    {
        /* start and end point are always a start and end point in one of the new segments */
        t.push_back(0.0) ;
        t.push_back(1.0) ;
        /* Sort the parameters numerically in ascending order 
        */
        sort( t.begin(), t.end()) ;

        /* Eliminate points that are the same, so are neighbors in the sorted list.
        * Work with the local copy: move upwards in the list and delete the points later in the list.
        */
        for( int refidx = 0 ; refidx < t.size()-1 ; )
        {
            if ( fabs(t[refidx]-t[refidx+1]) < pointFuzzy )
            {
                if ( refidx == 0 )
                    t.erase( t.begin() + refidx+1) ;
                else
                    t.erase(t.begin()+refidx) ;
            }
            else
                refidx ++ ;
        }

        /* Duplicates (pairs, triples,...) are now removed and pairs of adjacent parameters in the list
        * define the new line segments to be returned.
        */
        vector<Line> out ;
        for( int refidx = 0 ; refidx < t.size()-1 ; refidx++)
        {
            const Line seg(atParam( t[refidx]), atParam(t[refidx+1])) ;
            out.push_back(seg) ;
        }
        return out ;
    }

    /** Decide whether the parameter t that runs from 0 at the start
    * of the line segment up to 1 at the end is in the range, allowing
    * for some jitter in the floating point representation
    * @param t The line parameter.
    * @return true if the line parameter is inside the range (0,1).
    */
    static bool inSplitRange( const double t)
    {
        return ( t > pointFuzzy && t < 1.0-pointFuzzy ) ;
    }

    /** Decide whether the parameter t is close to one of the two values marking the terminal points.
    * @param t The line parameter.
    * @return true if the line parameter is close to 0 or 1.
    */
    static bool inEndRange( const double t)
    {
        return ( fabs(t) < pointFuzzy || fabs(t-1.0) < pointFuzzy ) ;
    }

protected:
private:
    static const double pointFuzzy = 1.e-10 ;

    /** Determinant of the four values within a 2 by 2 matrix.
    * @return the value of a00*a11 - a01*a10 .
    */
    static double determinant(const double a00, const double a01, const double a10, const double a11) 
    {
        return a00*a11-a01*a10 ;
    }
} ; /* Line */


/** A set of line segments
*/
class LineSet {
public:
    vector<Line> lines ;


    /** Intersect all lines in the set and split them such that 
    * all crossings are turned into terminal points of newly generated line elements
    * @return The new line set with no intersections left.
    */
    LineSet splitAtInters() const
    {
        /* The updated list of lines to be returned. The outbound set. */
        LineSet split;

        /* Take the existing line elements one by one and intersect them
        * with all elements of the outbound set.
        */
        for( int lidx =0 ; lidx < size() ; lidx++)
        {
            if ( split.size() == 0 )
                split += lines[lidx] ;
            else
            {
                /* The parameter list where lines[lidx] will be cut.
                */
                vector<double> tpoints ;

                /* The indices of elements in split[] that have been split by intersection with the current lidx .
                */
                set<int> rmPoints ;

                /* The subsegments of split[] that are created by cutting with the current lidx.
                */
                LineSet chaff ;

                /* loop over the line segments of the outbound set */
                for(int slidx = 0 ; slidx < split.size() ; slidx++)
                {
                    double t2 ;
                    /* do both lines intersect or touch ? */
                    double t = lines[lidx].inters(split.lines[slidx], t2) ;
                    if ( isinf(t) )
                        /* no action if parallel, since in the current problem these segments are distinct */
                        ;
                    else if ( Line::inSplitRange( t ) )
                    {
                        if ( Line::inSplitRange(t2) )
                        {
                            /* They do intersect: mark position t for splitting, mark element slidx
                            * for removal, and create its spliced segments for subsequent inclusion.
                            */
                            tpoints.push_back(t) ;
                            chaff += split.lines[slidx].splitAt(t2) ;
                            rmPoints.insert(slidx) ;
                        }
                        else if ( Line::inEndRange(t2) )
                            /* line in the outbound set touches the line at lidx: keep line in the
                            * outbound set and mark position t for subsequent splitting of lidx.
                            */
                            tpoints.push_back(t) ;
                    }
                    else if ( Line::inEndRange( t ) )
                    {
                        if ( Line::inSplitRange(t2) )
                        {
                            /* line lidx touches the line in the outbound set. Same replacement
                            * methodology for the outbound set as above, but keeping lidx intact.
                            */
                            chaff += split.lines[slidx].splitAt(t2) ;
                            rmPoints.insert(slidx) ;
                        }
                    }
                }

                /* Remove the outbound segments that have been split. Work in place which
                * requires working in reverse numerical order (from the maximum in rmPoints to the minimum)
                * to keep the indices intact.
                */
                while ( ! rmPoints.empty() )
                {
                    set<int>::iterator rmidx = max_element(rmPoints.begin(),rmPoints.end()) ;
                    /* remove the full line in the outbound set */
                    split.lines.erase(split.lines.begin() + *rmidx) ;

                    /* remove its index, so the next iteration of the loop will get the next smaller 
                    * index. */
                    rmPoints.erase(rmidx) ;
                }

                /* Insert the lines generated by splicing those removed above (list may be empty) */
                split += chaff ;

                /* Add the line lidx to the set, or its components if split markers had been created above.
                */
                if( tpoints.empty() )
                    split += lines[lidx] ;
                else
                    split += lines[lidx].splitAt(tpoints) ;
            }
        }

        /* return the updated list of segments, which has no crossings (only contacts) left */
        return split ;
    }

    /** Number of lines in the set.
    * @return number of edges in the current list.
    */
    int size() const
    {
        return lines.size() ;
    }

    /** Number of vertices in the set.
    * @return The number of vertices. This counts all terminal points, the duplicates
    * at vertices shared by one or more edges only counted once.
    */
    int vertices() const
    {
        set<Point> v;
        for(int lidx =0 ; lidx < lines.size() ; lidx++)
        {
            /* insertion of the two terminal points of this line.
            * The Point::operator== above eliminates duplicates at that time.
            */
            for(int i=0 ; i < 2 ; i++)
                v.insert( lines[lidx].pts[i]) ;
        }
        /* report the order of the set (count of element) finally residing in the point set 
        */
        return v.size() ;
    }

    /** Number of faces.
    * @return the count of tiles as computed by Euler's formula.
    */
    int faces() const
    {
        return 1+size()-vertices() ;
    }

    /** Add another line segment.
    * @param l The new additional line segment.
    */
    LineSet & operator += (const Line & l)
    {
        lines.push_back(l) ;
        return *this ;
    }

    /** Add another bundle of line segments.
    * @param l The list of additional line segments to be included.
    * @warn The algorithm does not check for duplications.
    */
    LineSet & operator += (const LineSet & l)
    {
        lines.insert(lines.end(), l.lines.begin(), l.lines.end() ) ;
        return *this ;
    }

    /** Add another bundle of line segments.
    * @param l The list of additional line segments to be included.
    * @warn The algorithm does not check for duplications.
    */
    LineSet & operator += (const vector<Line> & l)
    {
        lines.insert(lines.end(), l.begin(), l.end() ) ;
        return *this ;
    }

protected:
private:
} ; /* LineSet */

/** Regular polygon with an even number of edges.
*/
class Polyg : public LineSet {
public:
    /** Ctor with count of edges
    * @param n half of the number of edges.
    */
    Polyg(int n)
    {
        /* Construct the edges along the perimeter on the unit circle.
        */
        for(int e=0 ; e < 2*n; e++)
        {   
            Point a(cos((double)e*M_PI/n),sin((double)e*M_PI/n)) ;
            Point b(cos((1.+e)*M_PI/n),sin((1.+e)*M_PI/n)) ;
            *this += Line(a,b) ;
        }

        /* Add the diagonals parallel to sides.
        */
        for(int e=0 ; e < n; e++)
        {   
            int e2 = e+1 ;
            for(int k=1 ; k < n-1; k++)
            {
                Point a(cos((double)(e-k)*M_PI/n),sin((double)(e-k)*M_PI/n)) ;
                Point b(cos((double)(e2+k)*M_PI/n),sin((double)(e2+k)*M_PI/n)) ;
                *this += Line(a,b) ;
            }
        }
    }
protected:
private:

} ; /* Polyg */

/** Main program.
* @param argv There is one command line argument, the parameter n (half the number of edges)
*/
int main(int argc, char *argv[])
{
    int n = atoi(argv[1]) ;
    /* construct the polygon and the diagonals */
    Polyg p(n) ;
    /* Split all pairs of edges of this graph at all crossings */
    LineSet l = p.splitAtInters() ;
    /* report the number of edges, vertices and faces present */
    cout << l.size() << " edges " << l.vertices() << " vertices " << l.faces() << " tiles\n" ;

    return 0 ;
}
\end{verbatim}\normalsize

\bibliographystyle{amsplain}
\bibliography{all}

\end{document}